\documentclass[11pt]{amsart}

\usepackage{amsmath,amssymb,amsthm}
\usepackage{graphicx}

\def\dfrac{\displaystyle\frac}
\def\dsum{\displaystyle\sum}

\def\dlim{\displaystyle\lim}
\newtheorem{prop}{Proposition}
\newtheorem{theo}[prop]{Theorem}
\newtheorem{lemm}[prop]{Lemma}
\newtheorem{coro}[prop]{Corollary}
\newtheorem{rema}[prop]{Remark}

\newcommand{\ra}{\rightarrow}

\newcommand{\pa}{\partial}
\newcommand{\al}{\alpha}
\renewcommand{\leq}{\leqslant}
\renewcommand{\geq}{\geqslant}

\newcommand{\ep}{\varepsilon}

\numberwithin{equation}{section}

\title{A general form of Newton-Maclaurin type inequalities}

\begin{document}

\author{Changyu Ren}
\address{School of Mathematical Science\\
Jilin University\\ Changchun\\ China}
\email{rency@jlu.edu.cn}
\begin{abstract}
In this paper, we extend the classical Newton-Maclaurin inequalities
to functions $S_{k;s}(x)=E_k(x)+\dsum_{i=1}^s \al_i E_{k-i}(x)$,
which are formed by linear combinations of multiple basic symmetric
mean. We prove that for $\al_1,\al_2,\cdots,\al_s$ such that the
polynomial
$$t^s+\al_1 t^{s-1}+\al_2 t^{s-2}+\cdots+\al_s $$ has only real roots,
the Newton-Maclaurin type inequalities hold for $S_{k;s}(x)$.
\end{abstract}
\maketitle

\section{introduction}

The $k$-th elementary symmetric function of the variables
$x_1,x_2,\cdots,x_n$ is defined by
$$
\sigma_k(x)=\displaystyle\sum_{1\leq i_1<\cdots<i_k\leq
n}x_{i_1}\cdots x_{i_k}, \quad 1\leq k\leq n,
$$
where $x=(x_1,x_2,\cdots,x_n)$. It will be convenient to define
$\sigma_0(x)=1$, and define $\sigma_k(x)=0$ if $k<0$ or $k>n$.
Furthermore, define a $k$-th elementary symmetric mean as
$$
E_k(x)= \dfrac{\sigma_k(x)}{C_n^k},\quad k=0,1,\cdots,n,
$$
where $C_n^k=\dfrac{n!}{k!(n-k)!}$.

\par
The classic Newton inequalities \cite{Newt} and Maclaurin
inequalities \cite{Mac} are
\begin{equation}\label{e1.1}
E_k^2(x)\geq E_{k-1}(x)E_{k+1}(x),\quad k=1,2,\cdots,n
\end{equation}
and
\begin{equation}\label{e1.2}
E_1(x)\geq E_2^{1/2}(x)\geq\cdots\geq E_k^{1/k}(x),\quad 2\leq k\leq
n.
\end{equation}
Inequalities \eqref{e1.1} holds for any $x\in \mathbb{R}^n$, and
inequalities \eqref{e1.2} holds under the condition that $E_i\geq 0$
for all $i=1,\cdots, k$.

\par
A straightforward corollary of Newton's inequality \eqref{e1.1} is
\begin{equation}\label{e1.3}
\sigma_k^2(x)-\sigma_{k-1}(x)\sigma_{k+1}(x)\geq \theta
\sigma_k^2(x),\quad k=1,2\cdots,n,
\end{equation}
 where $\theta=\dfrac{(C_n^k)^2-C_n^{k-1}C_n^{k+1}}{(C_n^k)^2}$ is
a positive constant.

\par
The Newton-Maclaurin inequalities and its corollaries play important
roles in deriving theoretical result for fully nonlinear partial
differential equations and geometric analysis. There are many
important results that need to use the Newton-Maclaurin
inequalities, such as \cite{BCG,CX,BP,PX,PXF,GQ,GRW,GW,GC}, etc.
This is because the following $k$-Hessian equations and curvature
equations
\begin{align*}
\sigma_k(\lambda(u_{ij}))&=f(x,u,\nabla u), \\
\sigma_k(\kappa(X))&=f(X,\nu)
\end{align*}
are central topics in the field of fully nonlinear partial
differential equations and geometric analysis. The $k$-Hessian
operator $\sigma_k$ in the left-hand side  is the primary symmetric
function about the eigenvalues
$\lambda=(\lambda_1,\cdots,\lambda_n)$ of the Hessian matrix
$(u_{ij})$ or
 the principal curvature
$\kappa=(\kappa_1,\cdots,\kappa_n)$ of the surface.
\par

\par
 In recent years,
the fully non-linear differential equations derived from linear
combinations of primary symmetric functions have received increasing
attention. A significant example is the following special Lagrangian
equation
\begin{equation}\label{e1.4}
{\rm Im det}(\delta_{ij}+{\rm
i}u_{ij})=\displaystyle\sum_{k=0}^{[(n-1)/2]}(-1)^k\sigma_{2k+1}(\lambda(u_{ij}))=0,
\end{equation}
which was derived by Harvey and Lawson in their study of the minimal
submanifold problem \cite{HL}. They showed that if $u$ is a
solution, then the graph of $\nabla u$ is an absolutely
volume-minimizing submanifold of $\mathbb{R}^{2n}$. In addition,
there are many research works on nonlinear partial differential
equations where the operator on the left-hand side of the equation
is formed through linear combinations of multiple basic symmetric
functions, as seen in references \cite{Dong,Dws,GZ,Kry,LRW2,LR,Zhou}
and others.

\par
In the further study of the above problem, the Newton-Maclaurin type
inequalities for the functions derived from the left-hand side of
the equations are always needed.
 A natural question is whether the
Newton-Maclaurin type inequalities still hold for the functions of
linear combinations of these primary symmetric functions.

\par
In this paper, we will establish Newton-Maclaurin type inequalities
for the following functions
\begin{equation}\label{e1.5}
S_{k;s}(x)=E_k(x)+\dsum_{i=1}^s \al_i E_{k-i}(x), \quad 1\leq s<k<n,
\end{equation}
\begin{equation}\label{e1.6}
Q_{k;s}(x)=\sigma_k(x)+\dsum_{i=1}^s \al_i \sigma_{k-i}(x), \quad
1\leq s<k<n,
\end{equation}
where $x\in \mathbb{R}^n, \al=(\al_1,\cdots,\al_s)\in \mathbb{R}^s$,
and $S_{k;s}(x)$ and $Q_{k;s}(x)$ are, respectively, linear
combinations of multiple elementary symmetric means and elementary
symmetric functions. Specifically, we will establish inequalities of
the form \eqref{e1.1} and \eqref{e1.2} for $S_{k;s}(x)$, and
inequalities of the form \eqref{e1.3} for $Q_{k;s}(x)$. It is
noteworthy that the definitions of $S_{k;s}(x)$ and $Q_{k;s}(x)$
remain valid even when $s\geq k$. For example, when $s>2$, we have
\begin{align*}
&S_{1;s}(x)=E_1(x)+\al_1,\quad S_{2;s}(x)=E_2(x)+\al_1E_{1}+\al_2, \quad \cdots,\\
& S_{s;s}(x)=E_s(x)+\dsum_{i=1}^s\al_iE_{s-i}(x).
\end{align*}

\par
 In the recent work \cite{Ren} and \cite{HRW}, we have
established Newton-Maclaurin type inequalities for the cases $s=1$
and $s=2$ respectively. There are counterexamples in \cite{Ren} and
\cite{HRW} showing that inequalities \eqref{e1.1} does not always
hold  for $s>1$. Therefore, we need to impose some structural
conditions on $\al$.
\par
{\bf Condition C.} We say that $\al$ satisfies condition C if the
following $s$-degree polynomial
$$
f(t)=t^s+\al_1 t^{s-1}+\al_2 t^{s-2}+\cdots+\al_s
$$
related to $\al$ has only real roots. We denote these real roots as
$$-\beta=(-\beta_1,-\beta_2,\cdots,-\beta_s).$$

\par
The main result of this paper is stated as follows.

\begin{theo}\label{th1}
 For any $x\in
\mathbb{R}^n$, $\alpha\in \mathbb{R}^s$ and $1\leq s<n-1$, if
$\alpha$ satisfies condition C, then
\begin{equation}\label{e1.7}
S_{k;s}^2(x)\geq S_{k-1;s}(x)S_{k+1;s}(x),\quad k=s+1,\cdots,n-1.
\end{equation}
The inequalities are strict unless $n$ of the elements among
$x_1,x_2,\cdots,x_n,-\beta_1,-\beta_2$, $\cdots$, $-\beta_s$ are
equal or both sides of the inequalities are zero values.
\end{theo}

\par
\begin{rema}
When $s=2$ and $k=3$, a counterexample given in \cite{HRW} shows
that inequality \eqref{e1.7} does not always hold if $\al$ does not
satisfy condition C. It implies that in the case $s=2, k=3$,
condition C is a necessary and sufficient condition for the
inequality \eqref{e1.7} in Theorem \ref{th1} to hold.
\end{rema}

\begin{rema}
It can be verified that the coefficients of every $E_k(x)$ in the
special Lagrangian equation satisfy Condition C. In fact, according
to \cite{Yuan}, the coefficients of $E_k(x)$ in the special
Lagrangian equation correspond to the polynomial
\begin{align*}
\displaystyle\sum_{k=0}^{[(n-1)/2]}(-1)^k\sigma_{2k+1}(\frac{1}{t},\cdots,\frac{1}{t})=
n\arctan(\frac{1}{t})=0.
\end{align*}
The real roots of this polynomial are
$$
\cot\dfrac{0\pi}{2n},\pm\cot\dfrac{2\pi}{2n},\pm\cot
\dfrac{4\pi}{2n},\cdots,\pm\cot\dfrac{(n-1)\pi}{2n}
$$
when $n$ is odd,  and
$$
\pm \cot\dfrac{\pi}{2n}, \pm \cot\dfrac{3\pi}{2n},\pm
\cot\dfrac{5\pi}{2n}, \cdots, \pm\cot \dfrac{(n-1)\pi}{2n}
$$
when $n$ is even. Thus, the function in the special Lagrangian
equation satisfies inequality \eqref{e1.7}.
\end{rema}

The following corollary can be directly obtained from Theorem
\ref{th1}.

\begin{coro}
For $\alpha\in \mathbb{R}^s$ and any $k=s+1,\cdots,n-1$, if there
exists $x_0\in \mathbb{R}^n$ such that
\begin{equation*}
S_{k;s}^2(x_0)< S_{k-1;s}(x_0)S_{k+1;s}(x_0),
\end{equation*}
then the polynomial $f(t)$ has complex roots.
\end{coro}

\par
Similar as \cite{GLP} for $E_{k}(x)$, by \eqref{e1.7} we have
\begin{coro}\label{cor2}
With the hypothetical conditions as in Theorem \ref{th1} and
$s<l<k\leq n$, if
 \begin{equation}\nonumber
 S_{q;s}(x)\geq 0,\quad\text{for all}~q=l,\cdots,k-1,
  \end{equation}
then
\begin{align*}
&S_{l;s}(x)S_{k-1;s}(x) \geq S_{l-1;s}(x)S_{k;s}(x).
\end{align*}
\end{coro}

\par
Similar to the proof of the Maclaurin type inequalities for
$S_{k;2}(x)$ in \cite{HRW}, using inequalities \eqref{e1.7}, we can
obtain the following Maclaurin type inequalities for $S_{k;s}(x)$.
\par
\begin{theo}\label{th3}
The hypothetical conditions are the same as in Theorem \ref{th1}. If
we further assume that  $\beta_1\geq 0,\beta_2\geq 0, \cdots,
\beta_s \geq 0, E_1(x)\geq 0,E_2(x)\geq 0, \cdots, E_s(x)\geq 0$,
and
  \begin{equation*}
 S_{m;s}(x)\geq 0,\quad {\rm for~~ all}~~ m=s, s+1\cdots,k,
\end{equation*}
 then
\begin{equation}\label{e1.8}
S_{1;s}(x) \geq S_{2;s}^{1/2}(x) \geq \cdots\geq
S_{k;s}^{1/k}(x),\quad k=2,3,\cdots,n.
\end{equation}
\end{theo}

\par
For $Q_{k;s}(x)$, we have a result similar to inequalities
\eqref{e1.3}.

\begin{theo}\label{th4}
 For any $x\in
\mathbb{R}^n$, $\alpha\in \mathbb{R}^s$ and $1\leq s<n-1$. If
$\alpha$ satisfies condition C, then
\begin{equation}\label{e1.9}
Q_{k;s}^2(x)- Q_{k-1;s}(x)Q_{k+1;s}(x)\geq \theta Q_{k;s}^2(x),\quad
k=1,\cdots,n,
\end{equation}
where
$$
\theta=\dfrac{(C_{n+s}^k)^2-C_{n+s}^{k-1}C_{n+s}^{k+1}}{(C_{n+s}^k)^2}
$$
is a positive constant.
\end{theo}

\begin{rema}
If $\al$ does not satisfy Condition C, inequality \eqref{e1.9} is
not always true. For example, let
$Q_{3;2}(x)=\sigma_3(x)+\sigma_1(x)$. At this time, $\al=(0,1)$, and
the corresponding polynomial $t^2+1=0$ has no real roots. Choose
$x=(\frac{1}{3},\frac{1}{3},2,3)$, then we have
\begin{equation*}
Q_{3;2}^2(x)- Q_{2;2}(x)Q_{4;2}(x)=-\dfrac{10}{9}<0.
\end{equation*}
\end{rema}

\par
Similar to \cite{GLP} for $E_k(x)$, by \eqref{e1.9} we have the
following Corollary.

\begin{coro}\label{cor4}
The hypothetical conditions are the same as in Theorem \ref{th4}.
Let $s<l<k\leq n$, if
 \begin{equation}\nonumber
 Q_{q;s}(x)\geq 0,\quad\text{for all}~q=l,\cdots,k-1,
  \end{equation}
then
\begin{align*}
&Q_{l;s}(x)Q_{k-1;s}(x) \geq (1+\theta) Q_{l-1;s}(x)Q_{k;s}(x),
\end{align*}
where $0<\theta<1$ is a constants depending only on $n, k$ and $s$.
\end{coro}

\par
The proof method of Theorem 1 in this paper is different from that
in \cite{Ren} and \cite{HRW}, and it mainly uses the distribution of
real roots of polynomials related to $x\in \mathbb{R}^n$. The
structure of the paper is as follows: In Section 2, we mainly
consider the distribution of real roots of polynomials related to
$x$. The proofs of the main results, Theorem \ref{th1} and Theorem
\ref{th4}, will be given in Section 3.

\section{preliminary}

In this section and thereafter, we consistently assume that $x\neq
{\bf 0}$, because when $x={\bf 0}$, these inequalities
\eqref{e1.7}-\eqref{e1.9} clearly hold.  For the given
$x=(x_1,\cdots,x_n)\in \mathbb{R}^n$, we may assume that $P(t)$ is
an $n$-degree polynomial with real roots $x_{1},x_{2},\cdots,x_{n}$.
Then $P(t)$ can be represent as
 \begin{equation}\label{e2.1}
 P(t)=\prod_{i=1}^{n}(t-x_{i})=E_{0}(x)t^{n}-C_{n}^{1} E_{1}(x) t^{n-1}+C_{n}^{2} E_{2}(x) t^{n-2}-\cdots+(-1)^{n} E_{n}(x).
 \end{equation}

 The following lemma is an useful tool to prove Newton's
inequalities from \cite{GLP,NCP,Syl}.
\begin{lemm}\label{lem2.1}
If
\begin{equation}\nonumber
F(x,y)=c_{0}x^{n}+c_{1}x^{n-1}y+\cdots+c_{n}y^{n}
\end{equation}
is a homogeneous function of the n-th degree in x and y which has
all its roots ${x}/{y}$ real, then the same is true for all
non-identical $0$ equations
\begin{equation}\nonumber
\frac{\pa^{i+j}F}{\pa x^{i}\pa y^{j}}=0,
\end{equation}
obtained from it by partial  differentiation with respect to $x$ and
$y$. Further, if  $Q$ is one of these equations, and it has a
multiple root $\gamma$, then $\gamma$ is also a root, of
multiplicity one higher, of the equation from which $Q$ is derived
by differentiation.
\end{lemm}

We will apply Lemma \ref{lem2.1} to the homogeneous polynomial
related to $P(t)$ below,
\begin{equation}\nonumber
F(t,s)=E_{0}( x)t^{n}-C_{n}^{1} E_{1}( x) t^{n-1} s+C_{n}^{2} E_{2}(
x) t^{n-2} s^{2}-\cdots+(-1)^{n} E_{n}( x) s^{n}.
\end{equation}
Consider the derivative of $\dfrac{\partial F}{\partial t}$ and let
$s=1$, we obtain the following polynomial
 \begin{equation}\label{e2.2}
 P_1(t)=E_{0}(x)t^{n-1}-C_{n-1}^{1} E_{1}( x) t^{n-2}+C_{n-1}^{2} E_{2}( x) t^{n-3}-\cdots+(-1)^{n-1}
 E_{n-1}(x).
 \end{equation}
 Note that
$E_0(x)=1$, and by Lemma \ref{lem2.1}, the polynomial $P_1(t)$ has
$n-1$ real roots. Similarly, differentiate $F(t,s)$ with respect to
$s$ and let $s=1$, we obtain the following polynomial
 \begin{equation}\label{e2.3}
 P_2(t)=E_{1}(x)t^{n-1}-C_{n-1}^{1} E_{2}( x) t^{n-2}+C_{n-1}^{2} E_{3}( x) t^{n-3}-\cdots+(-1)^{n-1}
 E_{n}(x).
 \end{equation}
By Lemma \ref{lem2.1}, when $E_1(x)\neq 0$, the polynomial $P_2(t)$
has $n-1$ real roots.  When $E_{1}(x)=0$, using the identity
$$
[nE_{1}(x)]^2=\sigma_{1}^2(x)=\dsum_{i=1}^n
x_i^2+2\sigma_2(x)=\dsum_{i=1}^n x_i^2+2C_n^2 E_2(x),
$$
 it is easy to see that $E_{2}(x)\neq 0$. So in this case, the polynomial
$P_2(t)$ has $n-2$ real roots.

\par
Since
\begin{align*}
C_{n-1}^k+C_{n-1}^{k-1}=C_n^k,
\end{align*}
we have
\begin{equation}\label{e2.4}
P(t)=P_1(t)t-P_2(t).
\end{equation}

In the following, we will discuss the distribution of real roots of
the polynomial $P_2(t)$. Without loss of generality, we will always
assume $x_1\leq x_2\leq\cdots\leq x_n$ from now on.

\begin{lemm}\label{lem2.2}
Assume $x\in \mathbb{R}^n$ and $x_1< x_2<\cdots< x_n$, then the real
roots of polynomials $P_1(t)$ and $P_2(t)$ are all simple, and the
real roots of  $P_1(t)$ and  $P_2(t)$ are interleaved.
\end{lemm}

\begin{proof}
Since $x_1<x_2<\cdots<x_n$, that is, all roots of the polynomial
$P(t)$ are simple, by Rolle's Theorem, it is known that the real
roots $y_1,y_2,\cdots,y_{n-1}$ of $P_1(t)$ are located between each
pair of $x_1,x_2,\cdots,x_n$, as seen in Figure 1.

\begin{center}
$\begin{array}{c}
\includegraphics[scale=0.8]{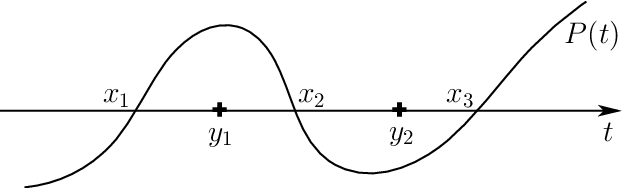}\\
\hbox{Figure 1}
\end{array}
$
\end{center}

We may consider a segment of the polynomial $P(t)$ curve for
discussion, which includes three adjacent real roots
$x_1',x_2',x_3'$ of $P(t)$ and two roots $y_1',y_2'$ of $P_1(t)$, as
shown in Figure 2.

\begin{center}
$\begin{array}{c}
\includegraphics[scale=0.8]{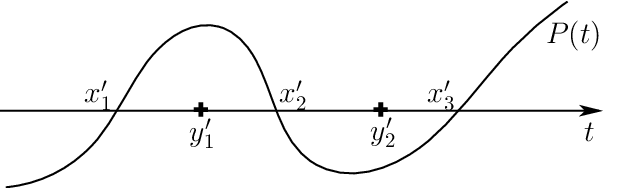}\\
\hbox{Figure 2}
\end{array}
$
\end{center}
\par
By \eqref{e2.4} we have
$$
P_2(y_1')=-P(y_1'),\quad P_2(y_2')=-P(y_2').
$$
It is evident that $P_2(t)$ has opposite signs at the points $y_1'$
and $y_2'$, therefore, there must be a real root of $P_2(t)$ between
$y_1'$ and $y_2'$. Similarly, we know that between each pair of
adjacent points $y_1,y_2,\cdots,y_{n-1}$, there is at least one real
root of $P_2(t)$. These roots add up to at least $n-2$ in total.

\par
We now prove that between any two adjacent real roots of $P_1(t)$,
there is at most one real root of $P_2(t)$. By contradiction, assume
that between two adjacent real roots $y_1'$ and $y_2'$ of $P_1(t)$,
there are two real roots $z_1$ and $z_2$ of $P_2(t)$. Let's assume
that the graph of $P_2(t)$ is shown as a dotted line in Figure 3.

\begin{center}
$\begin{array}{c}
\includegraphics[scale=0.8]{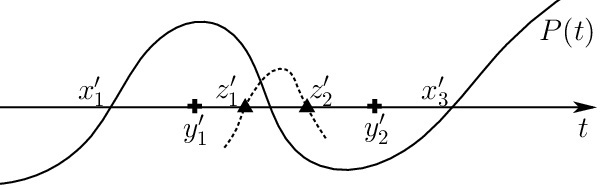}\\
\hbox{Figure 3}
\end{array}
$
\end{center}
From \eqref{e2.4}, we have $P_2(y_2')=-P(y_2')>0$. Note that $z_2'$
is a simple root of the polynomial $P_2(t)$, hence, there exists
some point $z'\in (z_2',y_2')$ such that $P_2(z')<0$. Therefore,
there must exist some point $z_0'\in(z',y_2')$ such that
$P_2(z_0')=0$. This means that $P_2(t)$ has at least three real
roots $z_1',z_2',z_0'$ in the interval $(y_1',y_2')$, which
contradicts the fact that $P_2(t)$ have at most $n-1$ real roots.
So, there is exactly one real root of $P_2(t)$ between each pair of
adjacent points $y_1'$ and $y_2'$, which indicates that the real
roots of polynomials $P_1(t)$ and $P_2(t)$ are interleaved.
\end{proof}

\par
For any real number $b$, we consider the polynomial
\begin{align}\label{e2.5}
P_3(t)=&P_2(t)+b P_1(t)\\
=&[E_{1}(x)+b E_{0}(x)]t^{n-1}-C_{n-1}^{1}[E_{2}(x)+b E_{1}(x)]
t^{n-2}\nonumber\\
&+C_{n-1}^{2}[E_{3}(x)+b E_{2}(x)] t^{n-3}-\cdots+(-1)^{n-1}
 [E_{n}(x)+b E_{n-1}(x)]. \nonumber
\end{align}

\begin{lemm}\label{lem2.3}
$\forall b\in \mathbb{R}, \forall x\in \mathbb{R}^n$, the polynomial
$P_3(t)$ has only real roots.
\end{lemm}
\begin{proof} When $ b=0$, $P_3(t)=P_2(t)$, and the conclusion is
obvious. Therefore, we assume that $ b\neq 0$. First, let's assume
$x_1<x_2<\cdots<x_n$ and discuss the cases where $E_1(x)>0$,
$E_1(x)<0$ and $E_1(x)=0$.

\par
{\bf Case A:} $E_1(x)>0$. In this case, let's further divide it into
the following sub-cases for discussion:

\par
{\bf Subcase A1:} $ b>0$ and $z_{n-1}>y_{n-1}$.
\par
Since $E_1(x)>0, b>0$, and when $t\rightarrow+\infty$,
$P_2(t)\rightarrow+\infty, P_1(t)\rightarrow+\infty$. The graphs of
the polynomials $ b P_1(x)$ and $P_2(x)$ are shown in Figure 4.

\begin{center}
$\begin{array}{c}
\includegraphics[scale=0.8]{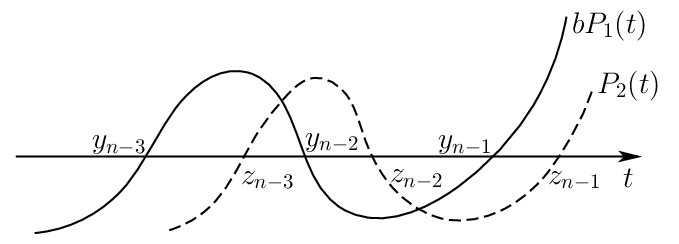}\\
\hbox{Figure 4}
\end{array}
$
\end{center}
\par\noindent
Since $P_3(y_{n-1})=P_2(y_{n-1})<0$ and $P_3(z_{n-1})= b
P_1(z_{n-1})>0$, there exists an $s_{n-1}\in (y_{n-1}, z_{n-1})$
such that $P_3(s_{n-1})=0$. Similarly, we can obtain that for each
$i=1,\cdots,n-2$, there exists an $s_i\in(y_i,z_i)$ such that
$P_3(s_{i})=0$, which meaning that $P_3(t)$ has only real roots.

\par
{\bf Subcase A2:} $ b>0$ and $z_{n-1}<y_{n-1}$.
\par
The proof is similar to Subcase A1; for each $i=1,\cdots,n-1$, there
exists an $s_i\in(z_i,y_i)$ such that $P_3(s_{i})=0$.

\par
{\bf Subcase A3:} $ b<0$ and $z_{n-1}>y_{n-1}$.

\par
Since $E_1(x)>0,  b<0$, when $t\rightarrow+\infty$,
$P_2(t)\rightarrow +\infty,  b P_1(t)\rightarrow -\infty$, thus the
graphs of $ b P_1(x)$ and $P_2(x)$ are as shown in Figure 5.

\begin{center}
$\begin{array}{c}
\includegraphics[scale=0.8]{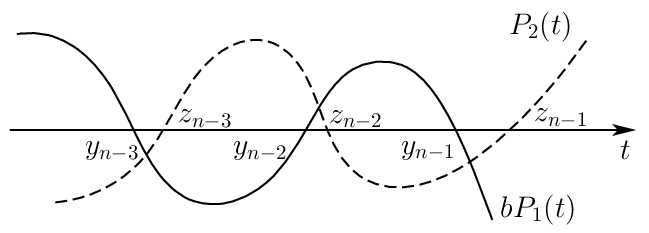}\\
\hbox{Figure 5}
\end{array}
$
\end{center}
\par\noindent
Since $P_3(z_{n-2})= b P_1(z_{n-2})>0$ and $P_3(y_{n-1})=
P_2(y_{n-1})<0$, there exists an $s_{n-2}\in (z_{n-2}, y_{n-1})$
such that $P_3(s_{n-2})=0$. Similarly, we can obtain that for each
$i=1,\cdots,n-3$, there exists an $s_i\in(z_i,y_{i+1})$ such that
$P_3(s_{i})=0$, which meaning that $P_3(t)$ has at least $n-2$ real
roots. Since $P_3(t)$ is an $n-1$-degree polynomial, $P_3(t)$ only
has real roots.

\par
{\bf Subcase A4:} $ b<0$ and $z_{n-1}<y_{n-1}$.

\par
The proof this subcase is similar to Subcase A3.

\par
{\bf Case B:} $E_1(x)<0$.
\par
The proof is similar to the Case A.

\par
{\bf Case C:} $E_1(x)=0$.
\par
In this case, the polynomial $P_2(t)$ has exactly $n-2$ real roots,
which are distributed between the real roots of $P_1(t)$, as stated
in Lemma \ref{lem2.2}. Let's assume $E_2(x)>0$ and discuss the two
cases where $ b>0$ and $ b<0$. The case where $E_2(x)<0$ is similar
to $E_2(x)>0$.

\par
{\bf Subcase C1:} $ b>0$.

\par
In this case, the graphs of polynomials $ b P_1(x)$ and $P_2(x)$ are
as shown in Figure 6.

\begin{center}
$\begin{array}{c}
\includegraphics[scale=0.8]{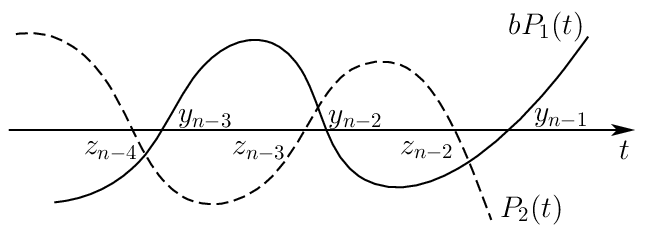}\\
\hbox{Figure 6}
\end{array}
$
\end{center}
\par\noindent
Since $P_3(z_{n-2})= b P_1(z_{n-2})<0$ and $P_3(y_{n-2})=
P_2(y_{n-2})>0$, there exists an $s_{n-2}\in (y_{n-2}, z_{n-2})$
such that $P_3(s_{n-2})=0$. Similarly, we can obtain that for each
$i=1,\cdots,n-3$, there exists an $s_i\in(y_i,z_{i})$ such that
$P_3(s_{i})=0$, which meaning that $P_3(t)$ has at least $n-2$ real
roots. Since $P_3(t)$ is an $n-1$-degree polynomial, $P_3(t)$ only
has real roots.

\par
{\bf Subcase C2:} $ b<0$.

\par
In this case, the graphs of polynomials $b P_1(x)$ and $P_2(x)$ are
as shown in Figure 7.
\begin{center}
$\begin{array}{c}
\includegraphics[scale=0.8]{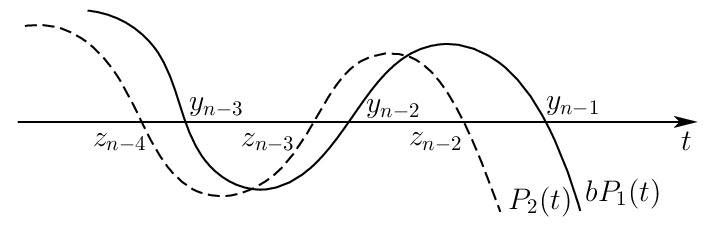}\\
\hbox{Figure 7}
\end{array}
$
\end{center}
\par\noindent
Similar to the previous discussion, for each $i=1,\cdots,n-2$, there
exists an $s_i\in(z_i,y_{i+1})$ such that $P_3(s_{i})=0$, which
means $P_3(t)$ has at least $n-2$ real roots. Since $P_3(t)$ is an
$n-1$-degree polynomial, $P_3(t)$ only has real roots.

\par
In summary, we proved that when $x_1<x_2<\cdots<x_n$, the polynomial
$P_3(t)$ has only real roots. Next, let's discuss the case where
$x_1\leq x_2\leq \cdots\leq x_n$.

\par
We assume that $x_1=x_2=\cdots=x_r<x_{r+1}$, for any
$r=1,2,\cdots,n-1$, then there exists a sufficiently small $\ep>0$
such that $x_1<x_2+\ep<x_3+2\ep<\cdots<x_r+(r-1)\ep<x_{r+1}$. We
denote $x^{\ep}_1=x_1, x^{\ep}_2=x_2+\ep, \cdots,
x^{\ep}_r=x_r+(r-1)\ep$. Using the same approach, for any
$x\in\mathbb{R}^n$, there exists a sufficiently small real number
$\ep>0$ such that $x^{\ep}_1<x^{\ep}_2<\cdots<x^{\ep}_n$. We denote
$x^\ep=(x^{\ep}_1, x^{\ep}_2,\cdots, x^{\ep}_n)$, then
$\dlim_{\ep\ra 0}x^{\ep}=x$.
\par

For the above $x^\ep$, denote the polynomials in Equations
\eqref{e2.1}-\eqref{e2.5} as $P^{\ep}(t)$, $P^{\ep}_1(t)$,
$P^{\ep}_2(t)$, $P^{\ep}_3(t)$. Applying Lemma \ref{lem2.2} and the
proof preceding this lemma, we know that the roots of polynomials
$P^{\ep}_1(t)$ and $P^{\ep}_2(t)$ are interlaced, and the polynomial
$P^{\ep}_3(t)$ has only real roots. Let $x\ra x^{\ep}$, we have
$$
P^{\ep}(t)\ra P(t),\quad P^{\ep}_1(t)\ra P_1(t), \quad
P^{\ep}_2(t)\ra P_2(t), \quad P^{\ep}_3(t)\ra P_3(t).
$$
Since the roots of a polynomial are continuous functions of its
coefficients, we have proven that for any $x\in \mathbb{R}^n$, the
polynomial $P_3(t)$ has only real roots.

\par
In fact, by Lemma 1, when $x_0$ is an $r$ repeated roots of $P(t)$,
$x_0$ is also an $r-1$ repeated roots of both $P_1(t)$ and $P_2(t)$.
Furthermore, from Equation \eqref{e2.5}, we can see that $x_0$
  is also an $r-1$-fold root of $P_3(t)$.
\end{proof}

\section{Proof of main results}

\subsection{Proof of Theorem \ref{th1} }

 We prove Theorem 1 by mathematical induction for $s=2,\cdots,s$. First, we have the following result (Theorem 1 in \cite{HRW}).

\begin{lemm}\label{lem2.4}
For any $x\in \mathbb{R}^n, \al=(\al_1,\al_2)\in \mathbb{R}^2$ and
$n>3$. If $\al$ satisfies condition C, then
\begin{equation*}
S^2_{k;2}(x)\geq S_{k-1;2}(x)S_{k+1;2}(x),\quad k=3,\cdots,n-1.
\end{equation*}
The inequality is strict unless both sides of the inequality are
zero or $n$ elements among $x_1, x_2, \cdots, x_n,
-\beta_1,-\beta_2$ are equal.
\end{lemm}

\par
Lemma \ref{lem2.4} shows that Theorem \ref{th1} holds when $s=2$. We
assume that Theorem \ref{th1} holds for $s-1$, that is

\begin{lemm}\label{lem2.5}
For any $x\in \mathbb{R}^n, \al'=(\al_1',\cdots,\al_{s-1}')\in
\mathbb{R}^{s-1}$ and $n>s$. If $\al'$ satisfies condition C, then
\begin{equation*}
[S'_{k;s-1}(x)]^2\geq S'_{k-1;s-1}(x)S'_{k+1;s-1}(x),\quad
k=s,\cdots,n-1.
\end{equation*}
Here $S'_{k;s-1}(x)=E_k(x)+\dsum_{i=1}^{s-1}\al'_iE_{k-i}(x)$, the
inequality is strict unless both sides of the inequality are zero
value or $n$ elements among $x_1, x_2, \cdots, x_n,
-\beta_1,\cdots,-\beta_{s-1}$ are equal.
\end{lemm}

\par
Next, we use Lemma \ref{lem2.5} and Lemma \ref{lem2.3} to prove
Theorem \ref{th1}. According to Lemma \ref{lem2.3}, for any $x\in
\mathbb{R}^n$, the polynomial $P_3(t)$ has only real roots. When
$E_{1}(x)+b E_{0}(x)\neq 0$, $P_3(t)$ has $n-1$ real roots, denoted
by $y=(y_1,y_2,\cdots,y_{n-1})$. Thus, $P_3(t)$ can be expressed as
\begin{align*}
P_3(t)=&\prod_{i=1}^{n-1}(t-y_{i})\\
=&t^{n-1}-C_{n-1}^{1}E_1(y) t^{n-2}+C_{n-1}^{2}E_{2}(y)
t^{n-3}-\cdots+(-1)^{n-1}
 E_{n-1}(y). \nonumber
 \end{align*}
 Comparing with the expression \eqref{e2.5} of
$P_3(t)$, it is evident that
\begin{align}\label{e3.4}
E_1(y)=&\dfrac{E_{2}(x)+ b E_{1}(x)}{E_{1}(x)+ b E_{0}(x)}, \qquad
E_2(y)=\dfrac{E_{3}(x)+ b E_{2}(x)}{E_{1}(x)+ b
E_{0}(x)},\\
\cdots,&\qquad  E_{n-1}(y)=\dfrac{E_{n}(x)+ b E_{n-1}(x)}{E_{1}(x)+
b E_{0}(x)}.\nonumber
\end{align}
For the above $y=(y_1,y_2,\cdots,y_{n-1})$, it follows from Lemma
\ref{lem2.5} that
\begin{align}\label{e3.5}
[S'_{k;s-1}(y)]^2\geq S'_{k-1;s-1}(y)S'_{k+1;s-1}(y),\quad
k=s,\cdots,n-2.
\end{align}

\par
Applying \eqref{e3.4}, we have
\begin{align}\label{e3.6}
S'_{k;s-1}(y)=&E_k(y)+\dsum_{i=1}^{s-1}\al_i'E_{k-i}(y)\nonumber\\
=&\dfrac{E_{k+1}(x)+ b E_{k}(x)}{E_{1}(x)+ b
E_{0}(x)}+\dsum_{i=1}^{s-1}\al_i'\dfrac{E_{k-i+1}(x)+ b
E_{k-i}(x)}{E_{1}(x)+ b E_{0}(x)}\nonumber\\
=&\dfrac{1}{E_{1}(x)+ b
E_{0}(x)}\Big[E_{k+1}(x)+(b+\al_1')E_{k}(x)+(b\al_1'+\al_2')E_{k-1}(x)+\cdots\\
&+(b\al_{s-2}'+\al_{s-1}')E_{k-s+2}(x)+b\al_{s-1}'E_{k-s+1}(x)
\Big]. \nonumber
\end{align}
Since $\al'=(\al_1',\cdots,\al_{s-1}')$ satisfies the condition C,
assume the polynomial
$$t^{s-1}+\al_1' t^{s-2}+\al_2' t^{s-3}+\cdots+\al_{s-1}' =0$$  has $s-1$ real roots $-\beta'=(-\beta_1,\cdots,-\beta_{s-1})$.
Then
$$
\al_{i}'=\sigma_{i}(\beta'), \quad i=1,2,\cdots, s-1.
$$
\par
Denote $\beta=(\beta_1,\cdots,\beta_{s-1},b)$ and let
$$\al_{i}=\sigma_{i}(\beta),\quad  i=1,2,\cdots,s.$$
Then $\al=(\al_1,\cdots,\al_{s})$ satisfies the condition C, and
$$
\al_1=b+\al_1',~ \al_2=b\al_{1}'+\al_{2}',~ \cdots,~
\al_{s-1}=b\al_{s-2}'+\al_{s-1}', ~\al_s=b\al_{s-1}'.
$$
Then equation \eqref{e3.6} can be rewritten as
\begin{align*}
S'_{k;s-1}(y)=&\dfrac{1}{E_{1}(x)+ b
E_{0}(x)}\Big[E_{k+1}(x)+\dsum_{i=1}^s\al_iE_{k+1-i}(x)\Big]\\
:=&\dfrac{S_{k+1;s}(x)}{E_{1}(x)+ b E_{0}(x)}.
\end{align*}
Similarly, we have
\begin{align*}
S'_{k-1;s-1}(y)=&\dfrac{S_{k;s}(x)}{E_{1}(x)+ b E_{0}(x)},\qquad
S'_{k+1;s-1}(y)=\dfrac{S_{k+2;s}(x)}{E_{1}(x)+ b E_{0}(x)}.
\end{align*}
 Thus, from \eqref{e3.5} we get
\begin{align*}
S_{k+1;s}^2(x)\geq S_{k;s}(x)S_{k+2;s}(x),\quad k=s,\cdots, n-2,
\end{align*}
or
\begin{align}\label{e3.7}
S_{k;s}^2(x)\geq S_{k-1;s}(x)S_{k+1;s}(x),\quad k=s+1,\cdots, n-1.
\end{align}

\par
When $E_{1}(x)+b E_{0}(x)=0$, let $b_{\ep}=b+\ep$, then
$E_{1}(x)+b_{\ep}E_{0}(x)=\ep\neq 0$. According to Lemma
\ref{lem2.3}, the polynomial
$$
P_3(t)=P_2(t)+b_{\ep} P_1(t)
$$
has only real roots, which further implies the following
inequalities
\begin{equation*}
[S^{\ep}_{k;s}(x)]^2\geq S^{\ep}_{k-1;s}(x)S^{\ep}_{k+1;s}(x),\quad
k=s+1,\cdots,n-1,
\end{equation*}
where
$$
S^{\ep}_{k;s}(x)=E_{k}(x)+\dsum_{i=1}^s\al_i^{\ep}E_{k-i}(x)
$$
with
$$\al_i^{\ep}=\sigma_i(\beta_1,\cdots,\beta_{s-1},b_{\ep}),\quad
i=1,2,\cdots, s. $$ Let $\ep\ra 0$, and we can obtain the
inequalities \eqref{e3.7}.

\par
Finally, we discuss the cases where equality holds in the
inequalities \eqref{e3.7}. By Lemma \ref{lem2.5}, when there are
$n-1$ equal elements among $y_1,\cdots, y_{n-1},
-\beta_1,\cdots,-\beta_{s-1}$, the equality in inequality
\eqref{e3.5} holds. Based on Lemma \ref{lem2.1}, it is also
equivalent to the condition that when $n$ equal elements among
$x_1,\cdots, x_{n}, -\beta_1,\cdots,-\beta_{s-1}$, the equality in
inequality \eqref{e3.7} holds. For each $\beta_i, i=1,\cdots,{s-1}$,
swap the positions of $b$ and $\beta_i$, and repeat the proof
process above, we can still obtain inequalities \eqref{e3.7}.
Moreover, the condition for equality in inequality \eqref{e3.7} is
that there are $n$ elements among $x_1,\cdots, x_{n},
-\beta_1,\cdots,-\beta_{i-1},-b,-\beta_{i+1},\cdots,-\beta_{s-1}$
are equal. When $\beta_1=\beta_2=\cdots=\beta_{s-1}=b$, according to
Theorem 3 of \cite{HRW}, the condition for equality in inequality
\eqref{e3.7} is that there are $n-s$ elements among $x_1,\cdots,
x_{n}$ are equal. Based on the above discussion, it can be concluded
that the condition for equality in inequality \eqref{e3.7} is that
there are $n$ equal elements among $x_1,\cdots,
x_{n},-\beta_1,\cdots,-\beta_{s-1},-b$. Thus, we proved Theorem
\ref{th1}.

\par
From the proof process of Theorem 1, we can easily obtain the
following corollary.
\begin{coro}
Let $1\leq s<k<n$. For the following real coefficients polynomial
$$
f(t)=t^n+C_n^1E_1t^{n-1}+C_n^2E_2t^{n-2}+\cdots+C_n^nE_n,
$$
denote $S_{k;s}=E_k+\dsum_{i=1}^s\al_iE_{k-i}$. If there exists an
$\al\in \mathbb{R}^s$ satisfying condition C such that
\begin{equation*}
S_{k;s}^2< S_{k-1;s}S_{k+1;s},
\end{equation*}
then the polynomial $f(t)$ has complex roots.
\end{coro}

\subsection{Proof of the Theorem \ref{th4}}

 Given
$x=(x_1,\cdots,x_n)\in \mathbb{R}^n$ and
$\beta=(\beta_1,\cdots,\beta_s)\in \mathbb{R}^s$, consider the
following $n$-degree polynomial
 \begin{equation}\label{e3.8}
 P(t)=\prod_{i=1}^{n}(t-x_{i})=\sigma_{0}(x)t^{n}-\sigma_{1}(x) t^{n-1}+\sigma_{2}(x) t^{n-2}-\cdots+(-1)^{n} \sigma_{n}(x)
 \end{equation}
with real roots $x_1,x_2,\cdots,x_n$. Obviously, the coefficients
$\sigma_k(x)$ $(k=1,2,\cdots,n)$ of this polynomial satisfy
inequalities \eqref{e1.3}. Let
 \begin{align}\label{e3.9}
 Q_1(t)=&(t-\beta_1)\prod_{i=1}^{n}(t-x_{i})=tP(t)-\beta_1 P(t) \\
 =&\sigma_{0}(x)t^{n+1}-[\sigma_{1}(x)+\beta_1\sigma_{0}(x)] t^{n}+[\sigma_{2}(x)+\beta_1\sigma_{1}(x)] t^{n-1}-\cdots  \nonumber\\
 &+(-1)^{n} [\sigma_{n}(x)+\beta_1\sigma_{n-1}(x)]t+(-1)^{n+1}
 \beta_1\sigma_{n}(x).\nonumber
 \end{align}
The polynomial $Q_1(t)$ of degree $n+1$ obviously has $n+1$ real
roots: $\beta_1, x_1, x_2,\cdots, x_n$.

\par
Denote $Y_1=(\beta_1, x_1, x_2,\cdots, x_n)$. Then
\begin{align}\label{e3.10}
Q_1(t)=\sigma_{0}(Y_1)t^{n+1}-\sigma_{1}(Y_1) t^{n}+\sigma_{2}(Y_1)
t^{n-1}-\cdots+(-1)^{n+1} \sigma_{n+1}(Y_1).
\end{align}
By comparing the coefficients of the polynomial $Q_1(t)$ in
\eqref{e3.9} and \eqref{e3.10}, we have
$$
\sigma_k(Y_1)=\sigma_{k}(x)+\beta_1\sigma_{k-1}(x),\quad  0\leq
k\leq n+1.
$$
Next, let
\begin{align*}
 Q_2(t)=&(t-\beta_2)(t-\beta_1)\prod_{i=1}^{n}(t-x_{i})=tQ_1(t)-\beta_2
 Q_1(t),
\end{align*}
then $Q_2(t)$ is a polynomial of degree $n+2$, and it has $n+2$ real
roots $\beta_1,\beta_2,x_1, \cdots,x_n$.
\par
 Similarly, denote
$Y_2=(\beta_1,\beta_2, x_1, \cdots, x_n)$. Then we have
\begin{align*}
\sigma_k(Y_2)=&\sigma_{k}(Y_1)+\beta_1\sigma_{k-1}(Y_1)\\
=&\sigma_{k}(x)+(\beta_1+\beta_2)\sigma_{k-1}(x)+\beta_1\beta_2\sigma_{k-2}(x),\quad
0\leq k\leq n+2.
\end{align*}
Repeating the above process $s$ times, we obtain
$Y_s=(\beta_1,\cdots, \beta_s, x_1, \cdots, x_n)$
 and
\begin{align*}
\sigma_k(Y_s)=&\sigma_{k}(x)+\sigma_1(\beta)\sigma_{k-1}(x)+\sigma_2(\beta)\sigma_{k-2}(x)+\cdots+\sigma_s(\beta)\sigma_{k-s}(x)\\
=&\sigma_{k}(x)+\dsum_{i=1}^s\al_i\sigma_{k-i}(x)=Q_{k;s}(x),\quad
0\leq k\leq n+s.
\end{align*}
\par
Applying inequality \eqref{e1.3} to $Y_s$ above, we obtain
\begin{align*}
[Q_{k;s}(x)]^2\geq (1+\theta)Q_{k-1;s}(x)Q_{k+1;s}(x),\quad 1\leq k<
n+s,
\end{align*}
where
$$
\theta=\dfrac{(C_{n+s}^k)^2-C_{n+s}^{k-1}C_{n+s}^{k+1}}{(C_{n+s}^k)^2}.
$$
Thus, we proved Theorem \ref{th4}.

\bigskip

\noindent {\it Acknowledgement:} The author wish to thank  Professor
Pengfei Guan and Zhizhang Wang for their valuable suggestions and
comments. Part of the work was done while the authors were visiting
The Chinese University of Hong Kong. I would like to
 thank  their support and hospitality.

\end{document}